\def\timestamp{%
Time-stamp: <soft-compactifications-topproc.tex: Monday 10-05-2021 at 21:46:26 (cest)>}
\def\stripname Time-stamp: <#1 #2>{#2}
\edef\filedate{\expandafter\stripname\timestamp}

\documentclass{amsart}
\usepackage[T1]{fontenc}   

\newtheorem{theorem}{Theorem}

\theoremstyle{definition}

\newtheorem{remark}[theorem]{Remark}
\numberwithin{equation}{section}

\theoremstyle{definition}
\newtheorem{example}{Example}
\newtheorem{question}{Question}

\DeclareMathSymbol\N  0{AMSb}{`N}

\DeclareMathSymbol\le    \mathrel{AMSa}{"36}
\DeclareMathSymbol\ge    \mathrel{AMSa}{"3E}
\DeclareMathSymbol\restr \mathbin{AMSa}{"16}

\newcommand\bee{\mathfrak{b}}
\newcommand\cee{\mathfrak{c}}
\newcommand\pee{\mathfrak{p}}
\newcommand\tee{\mathfrak{t}}

\newcommand\cl{\operatorname{cl}}
\newcommand\orpr[2]{\langle#1,#2\rangle}
\newcommand\preim{^\gets}

\newcommand\axiom{\mathsf}
\newcommand\CH{\axiom{CH}}
\newcommand\ZFC{\axiom{ZFC}}

\newcommand\betaN{\beta\N}
\newcommand\betaNminN{\betaN\setminus\N}

\newcommand\calA{\mathcal{A}}
\newcommand\calB{\mathcal{B}}
\newcommand\calF{\mathcal{F}}
\newcommand\pow{\mathcal{P}}
\newcommand\calL{\mathcal{L}}
\newcommand\calR{\mathcal{R}}
\newcommand\calS{\mathcal{S}}
\newcommand\calU{\mathcal{U}}
\newcommand\calX{\mathcal{X}}

\newcommand\omegaseq[2][\omega]{\langle{#2}_n:n<#1\rangle}
\newcommand\omegaoneseq[2][\omega_1]{\langle{#2}_\alpha:\alpha<#1\rangle}
\newcommand\omegaoneseqorpr[2]
          {\bigl<\orpr{#1_\alpha}{#2_\alpha}:\alpha<\omega_1\bigr>}
\newcommand\Ccube[1][\omega_1]{2^{#1}}
\newcommand\Tcube[1][\omega_1]{[0,1]^{#1}}

\usepackage{amsrefs}
\begin{document}

\noindent                                             
\begin{picture}(150,36)                               
\put(5,20){\tiny{Submitted to}}                       
\put(5,7){\textbf{Topology Proceedings}}              
\put(0,0){\framebox(140,34){}}                        
\put(2,2){\framebox(136,30){}}                        
\end{picture}                                        
\vspace{0.5in}

\renewcommand{\bf}{\bfseries}
\renewcommand{\sc}{\scshape}
\vspace{0.5in}

\title[Soft-Parovichenko spaces]%
      {All Parovichenko spaces may be soft-Parovichenko}

\author[A. Dow]{Alan Dow}
\address{
Department of Mathematics\\
UNC-Charlotte\\
9201 University City Blvd. \\
Charlotte, NC 28223-0001}
\email{adow@uncc.edu}
\urladdr{https://webpages.uncc.edu/adow}

\author[K. P. Hart]{Klaas Pieter Hart}

\address{Faculty EEMCS\\TU Delft\\
         Postbus 5031\\2600~GA {} Delft\\the Netherlands}
\email{k.p.hart@tudelft.nl}
\urladdr{http://fa.its.tudelft.nl/\~{}hart}



\date{\filedate}

\begin{abstract}
It is shown that, assuming the Continuum Hypothesis,
every compact Hausdorff space of weight at most~$\cee$ is a remainder
in a soft compactification of~$\N$.

We also exhibit an example of a compact space of weight~$\aleph_1$ ---
hence a remainder in some compactification of~$\N$ ---
for which it is consistent that is not the remainder in a soft
compactification of~$\N$.
\end{abstract}

\dedicatory{To the memory to Phil Zenor, one of founders of this journal}

\maketitle

\section*{Introduction}

A compactification, $\gamma\N$, of the discrete space~$\N$ of natural numbers
is said to be \emph{soft} if for all pairs~$\orpr AB$ 
of disjoint subsets of~$\N$ the following holds:
if $\cl A\cap\cl B\neq\emptyset$
then there is an autohomeomorphism~$h$ of~$\gamma\N$ such that
$h[A]\cap B$ is infinite and $h$~is the identity 
on the remainder~$\gamma\N\setminus\N$.
 
Banakh asked in~\cite{flow309583} whether every Parovichenko
space is soft-Parovi\-chen\-ko,
where a Parovichenko space is defined to be a remainder in some 
compactification of~$\N$ and, naturally, a soft-Parovichenko
space is a remainder in some soft compactification of~$\N$.
Parovichenko's classic theorem, from~\cite{MR0150732},  
characterizes, assuming $\CH$, the
Parovichenko spaces as the compact Hausdorff spaces of weight at 
most~$\cee$.

\begin{example}
The \v{C}ech-Stone compactification, $\betaN$, of~$\N$ is soft, vacuously;
hence $\betaNminN$~is soft-Parovichenko.

At the other end of the spectrum the one-point compactification~$\alpha\N$ 
is soft too as \emph{every} permutation of~$\N$ determines an 
autohomeomorphism of~$\alpha\N$.
\end{example}

\begin{remark} \label{rem.FU}
As remarked by Banakh in~\cite{flow309583}: 
if $\delta\N$ is a compactification of~$\N$ with the property that
whenever $x\in\cl A$, 
where $x\in\delta\N\setminus\N$ and $A\subseteq\N$,
there is a sequence in~$A$ that converges to~$x$,
then $\delta\N$~is a soft compactification. 

Indeed, if $x\in\cl A\cap\cl B$ and $S$ and~$T$ are subsets of~$A$ 
and~$B$ respectively that converge to~$x$ then one takes a permutation~$h$ 
of~$\N$ that interchanges $S$ and $T$ and is the identity outside $S\cup T$.
The extension of~$h$ by the identity on~$X$ is an autohomeomorphism.
\end{remark}

In \cite{MR4142223} the reader can find more information
on the background of the problem of remainders in soft compactifications,
including, in Theorem~9.5, some classes of compact spaces that are 
soft-Parovichenko:
\begin{itemize}
\item Parovichenko and of character less than $\pee$
\item perfectly normal
\item of weight less than $\pee$
\end{itemize}
In each case one obtains the stronger statement that every compactification 
with the space as its remainder is soft, because in each case the 
compactification satisfies the condition in Remark~\ref{rem.FU}.
The cardinal number~$\pee$ is equal to the cardinal number~$\tee$
discussed in Section~\ref{sec.examples} below.

\section{Applying the Continuum Hypothesis}

In this section we prove the statement in the abstract.
The Continuum Hypothesis ($\CH$) implies that every Parovichenko
space is soft-Parovichenko.

\smallskip
Let $X$ be compact Hausdorff and of weight~$\aleph_1$.
We may assume $X$ is embedded in the Tychonoff cube $[0,1]^{\omega_1}$ and,
for technical convenience, that $X\subseteq\{0\}\times[0,1]^{[1,\omega_1)}$.

Our aim will be to construct a sequence 
$\omegaoneseq{f}$ of functions from~$\N$ to~$[0,1]$
such that the \v{C}ech-Stone extension~$\beta f$ of its diagonal 
map $f:\N\to[0,1]^{\omega_1}$ satisfies $\beta f[\betaNminN]=X$.
To make sure that $f[\N]$~is discrete we demand that $f_0(n)=2^{-n}$ 
for all~$n$.
In this way $f[\N]\cup X$ will be a compactification of~$\N$ with $X$ as 
its remainder.
 
\subsection*{Ensuring softness}

To ensure softness of the compactification we take our inspiration from
Remark~\ref{rem.FU}.

Along with the functions $f_\alpha$ we construct an almost disjoint 
family~$\calS$ of subsets of~$\N$ such that in the end every~$S\in\calS$
converges to a point~$x_S$ of~$X$.
In addition we ensure that whenever $\orpr AB$ is a pair of disjoint
subsets of~$\N$ whose closures intersect then there will two sets~$S$ and~$T$
in~$\calS$ such that $S\cap A$ and $T\cap B$ are infinite, and $x_S=x_T$.
As in Remark~\ref{rem.FU} a permutation of~$\N$ that interchanges
$S\cap A$ and $T\cap B$ and is the identity outside these sets
gives an autohomeomorphism of the compactification as required.

\subsection*{The construction}

We let $\omegaoneseqorpr{A}{B}$ enumerate the set of ordered pairs of 
disjoint infinite subsets of~$\N$.
We shall construct
\begin{itemize}
\item a sequence $\omegaoneseq{f}$ of functions from $\N$ to $[0,1]$
\item a sequence $\omegaoneseq{S}$ of subsets of~$\N$
\item a sequence $\omegaoneseq{x}$ of points in~$X$
\item a sequence $\omegaoneseq{K}$ of subsets of~$\N$
\end{itemize}
For each $\delta$ we let $g_\delta:\N\to[0,1]^\delta$ be the diagonal map
of $\omegaoneseq[\delta]{f}$, and, for bookkeeping purposes,
$I$~will be the set of~$\delta$ for which the closures 
of~$g_\delta[A_\delta]$ and~$g_\delta[B_\delta]$
intersect.

The sequences should satisfy the following conditions.
\begin{enumerate}
\item if $\alpha<\delta$ then the set $g_\delta[S_\alpha]$ converges to
      the point~$x_\alpha\restr\delta$
\item if $\delta\in I$ then there are $\alpha,\beta\le\delta$ such that
      $x_\alpha\restr\delta=x_\beta\restr\delta$, and\label{cond.intersect}
      both intersections $S_\alpha\cap A_\delta$ and $S_\beta\cap B_\delta$ 
      are infinite
\item for all $\delta$ the family $\{S_\alpha:\alpha<\delta\}\cup\{K_\delta\}$
      is almost disjoint
\item if $\alpha<\beta$ then $K_\beta\subseteq^* K_\alpha$
\item if $\alpha\in\omega_1$ 
      then $\beta g_\alpha[\N^*]=\beta g_\alpha[K_\alpha^*]=X\restr\alpha$.
\end{enumerate}
In condition~\ref{cond.intersect} we do not exclude the possibility
that $\alpha=\beta$.

At each stage $\delta$ we choose the set~$S_\delta$,
construct the function~$f_\delta$, 
and determine the set~$K_{\delta+1}$ as a subset of~$K_\delta$.
This means that in case $\delta$~is a limit we must construct~$K_\delta$ 
first.

\subsubsection*{Making $K_\delta$ if $\delta$ is a limit}

Let $\omegaseq\calU$ be a sequence
of finite families of basic open sets in~$[0,1]^\delta$ such that
for all~$n$ we have $X\restr\delta\subseteq\bigcup\calU_n$ 
and $(X\restr\delta)\cap U\neq\emptyset$ for all~$U\in\calU_n$,
and such that for every open set~$O$ around~$X\restr\delta$ there is an~$n$
such that $\bigcup\calU_n\subseteq O$.
For every $n$ there is a finite set~$F_n$ such that every member of~$\calU_n$
has its support in~$F_n$.
Let $\omegaseq\delta$ be a strictly increasing sequence
of ordinals that converges to~$\delta$ and such that $F_n\subseteq\delta_n$
for all~$n$.

The family $\calU_n$ can also be considered to be a family of basic open sets 
in the product $[0,1]^{\delta_n}$.
The condition that $\beta g_{\delta_n}[K_{\delta_n}^*]=X\restr\delta_n$ 
for all~$n$ translates into two things:
\begin{itemize}
\item for every~$n$ there is a natural number~$N_n$ such 
      that $g_{\delta_n}(k)\in\bigcup\calU_n$ for $k\in K_{\delta_n}\setminus N_n$
\item for every $U\in\calU_n$ the set $\{k\in K_{\delta_n}:g_{\delta_n}(k)\in U\}$
      is infinite
\end{itemize}
But then the same holds with $g_\delta$ replacing $g_{\delta_n}$.

Using this we determine a strictly increasing sequence 
$\omegaseq M$ of natural numbers 
such that $M_n\ge N_n$ for all~$n$, 
such that $K_{\delta_{n+1}}\setminus M_{n+1}\subset K_{\delta_n}$ for all~$n$, 
and such that for every $U\in\calU_n$ 
    there is a $k\in K_{\delta_n}\cap[M_n,M_{n+1})$
such that $g_\delta(k)\in U$.

We let $K_\delta=\bigcup_{n<\omega}\bigl(K_{\delta_n}\cap[M_n,M_{n+1})\bigr)$.
By construction $\beta g_\delta[K_\delta^*]=X$, and because 
$K_\delta\subseteq^*K_{\delta_n}$ for all~$n$ the set is also almost disjoint
from all~$S_\alpha$ for $\alpha<\delta$.

\subsubsection*{The actual construction}

Now let $\delta\in\omega_1$ and assume that everything has been constructed
up to and\slash or including~$\delta$.

If the closures of $g_\delta[A_\delta]$ and $g_\delta[B_\delta]$ intersect
then we add $\delta$ to the set~$I$ and 
determine $S_\delta$ by considering a few cases.

First shrink $A_\delta$ and $B_\delta$ to infinite sets $C$ and $D$ such that 
the closures of $g_\delta[C]$ and $g_\delta[D]$ intersect in exactly one point
of~$X\restr\delta$, this point is going to grow into~$x_\delta$, so we denote
it~$x_\delta\restr\delta$.
Note that the union $g_\delta[C]\cup g_\delta[D]$ converges 
to~$x_\delta\restr\delta$.

The cases that can occur are
\begin{itemize}
\item both $C$ and $D$ are almost disjoint from the $S_\alpha$ 
      with $\alpha<\delta$; in this case we let $S_\delta$ be an infinite 
      subset of $C\cup D$, that meets both $C$ and $D$ in an infinite set
      and is such that $K_\delta\setminus S_\delta$ contains an infinite set
      that converges to~$x_\delta\restr\delta$.
\item $C$ is almost disjoint from the $S_\alpha$ with $\alpha<\delta$, but
      $D$~is not; in this case we have a $\beta<\delta$ such 
      that $S_\beta\cap D$ is infinite, and so 
      $x_\beta\restr\delta=x_\delta\restr\delta$. Now let $S_\delta$ be an 
      infinite subset of~$C$ as in the previous case.
\item $D$ is almost disjoint from the $S_\alpha$ with $\alpha<\delta$, but
      $C$~is not; in this case we have an $\alpha<\delta$ such 
      that $S_\alpha\cap D$ is infinite, and so 
      $x_\alpha\restr\delta=x_\delta\restr\delta$. 
      Now let $S_\delta$ be an infinite subset of~$D$ as in the previous cases.
\item neither $C$ nor $D$ is almost disjoint from the $S_\alpha$ 
       with $\alpha<\delta$; this means that condition~\ref{cond.intersect}
       is already met. 
      We let $S_\delta$ be an infinite subset of~$K_\delta$ that converges
      to some $x_\delta\restr\delta$, again subject to the condition from
      the first three cases.
\end{itemize}
In all four cases let $K_{\delta+1}=K_\delta\setminus S_\delta$.
Then $\beta g_\delta[K_{\delta+1}]=X$, because of the condition on~$S_\delta$.
That condition is met automatically if $x_\delta\restr\delta$~is not an isolated
point of~$X\restr\delta$.

\smallbreak
In case the closures do not intersect we choose an arbitrary 
point~$x_\delta\restr\delta$ of~$X\restr\delta$ and choose an infinite
subset~$S_\delta$ of~$K_\delta$ that converges to this point
and such that $K_\delta\setminus S_\delta$, which will be~$K_{\delta+1}$,
contains an infinite set that converges to~$x_\delta\restr\delta$ as well.

In all cases we still have $\beta g_\delta[K_{\delta+1}^*]=X$.

\smallbreak

Before we proceed to the definition of~$f_\delta$ we first choose the 
$\delta$th coordinates of the points~$x_\alpha$ for $\alpha\le\delta$.
For each $\alpha$ we check whether there is a $\beta<\alpha$
such that $x_\alpha\restr\delta=x_\beta\restr\delta$.
If that is the case then we must let $x_\alpha(\delta)=x_\beta(\delta)$.
In the other case we ensure that
$\orpr{x_\alpha\restr\delta}{x_\alpha(\delta)}\in X\restr(\delta+1)$. 
This then introduces the demand that $f_\delta[S_\alpha]$ 
converge to~$x_\alpha(\delta)$.

\smallbreak 
To specify the function~$f_\delta$ we proceed much as in the construction
of~$K_\delta$ for limit~$\delta$.

We take a sequence $\omegaseq\calU$ of finite families 
of basic open sets in~$[0,1]^{\delta+1}$ as follows.
First take an increasing sequence $\omegaseq F$
of finite subsets of~$\delta+1$ such that $\delta\in F_0$
and $\bigcup_nF_n=\delta+1$.

Next we let $\calB_n$ be the family of all products
$\prod_{\alpha\le\delta}I_\alpha$, where $I_\alpha$ is an interval of the
form $[0,2^{-n})$, $\bigl(i\cdot2^{-n},(i+1)\cdot2^{-n}\bigr)$
or $(1-2^{-n},1]$ in~$[0,1]$ if $\alpha\in F_n$
and $I_\alpha=[0,1]$ if $\alpha\notin F_n$.

We let $\calU_n=\{B\in\calB_n:B\cap(X\restr(\delta+1))\neq\emptyset\}$
and we write every $U\in\calU_n$ as $V_U\times I_U$, where $V_u$~is
in $[0,1]^\delta$ and $I_U$ is an interval in~$[0,1]$.

For every~$n$ we let  
$\calS_n=\{S_\alpha:\alpha\in F_n\}\cup\{K_{\delta+1}\}$ and we take 
an $N_n$ such that 
\begin{itemize}
\item for all distinct $X$ and $Y$ in $\calS_n$ the 
      intersection $X\cap Y$ is contained in~$N_n$
\item for all $k\ge N_n$ we have $g_\delta(k)\in\bigcup\{V_U:U\in\calU_n\}$
\item for all $U\in\calU_n$ and all $\alpha\le\delta$: 
      if $x_\alpha\restr\delta \in V_U$ then 
      $g_\delta(k)\in V_U$ for all $k\in S_\alpha\setminus N_n$ 
\end{itemize}

Because $\beta g_\delta[K_{\delta+1}^*]=X$ we know that for every $U\in\calU_n$
the set $\{k\in K_{\delta+1}:g_\delta(k)\in V_U\}$ is infinite.
Using this we take a strictly increasing sequence $\omegaseq M$
of natural numbers such that $M_n\ge N_n$ for all~$n$
and we can define $f_\delta$ on~$K_{\delta+1}$ 
such that for all~$n$ and $U\in\calU_n$ 
there a $k\in K_{\delta+1}\cap[M_n,M_{n+1})$ such that 
$\orpr{g_\delta(k)}{f_\delta(k)}\in U$.

\smallbreak
We define $f_\delta$ on $S_\alpha\cap[M_n,M_{n+1})$ whenever
$\alpha\in F_n$; because $M_n\ge N_n$ there will be no interference
with the values that we specified on~$K_{\delta+1}$ and between 
the different $S_\alpha$s.

For each $\alpha\in F_n$ we can simply define $f_\delta(k)=x_\alpha(\delta)$ 
for $k\in S_\alpha\cap[M_n,M_{n+1})$.

For all $k\in [M_n,M_{n+1})$ that are not in $K_{\delta+1}$,
nor in in any of the $S_\alpha$ for some $\alpha\le\delta$,
we simply choose $f_\delta(k)$ in such a way that 
$\orpr{g_\delta(k)}{f_\delta(k)}\in\bigcup\calU_n$.

To see that $f_\delta[S_\alpha]$ converges to~$x_\alpha(\delta)$ it suffices to 
observe that $f_\delta$ has the constant value~$x_\alpha(\delta)$ on the 
intersection $S_\alpha\cap[M_n,\omega)$, where $n$~is such 
that $\alpha\in F_n$.

\section{Some examples of weight $\aleph_1$}
\label{sec.examples}

One half of Parovichenko's characterization is a $\ZFC$ result: 
every compact space of weight~$\aleph_1$ is a remainder of~$\N$.
Our proof is in the spirit of B\l aszczyk and Szymanski's proof
of this statement in~\cite{MR628044}.
The main difference is in the number of tasks that need to be done
to construct the map $f:\N\to[0,1]^{\omega_1}$.
To make $X$ a remainder involves just $\aleph_1$ many tasks, whereas
to make it a soft remainder seems to involve $\cee$ many tasks as
there are that many pairs of subsets of~$\N$,
hence the need to assume~$\CH$.

In this section we present some examples to show that without $\CH$
not every compact space of weight~$\aleph_1$ is automatically a soft remainder
of~$\N$.

\subsection*{The ordinal space $\omega_1+1$}

This relatively simple space already offers some difficulties regarding 
softness.
To be sure: it is a soft remainder, but the proof requires us to consider
two cases, depending on the value of the small cardinal~$\tee$.
\begin{itemize}
\item if $\tee>\aleph_1$ then \emph{every} compactification of~$\N$
with $\omega_1+1$ as a remainder is soft, 
\item if $\tee=\aleph_1$ then \emph{some but not all} compactifications 
of~$\N$ with $\omega_1+1$ as a remainder are soft.
\end{itemize}

All compactifications of~$\N$ with~$\omega_1+1$ as a remainder have roughly
the same structure, as described by Franklin and Rajagopalan in~\cite{MR283742}.

Let $\gamma\N$ be such a compactification, where we assume that $\N$ 
and $\omega_1+1$ are disjoint.
We apply normality to find, for every $\alpha<\omega_1$.
open sets $U_\alpha$ and $V_\alpha$ with disjoint closures such 
that $[0,\alpha]\subseteq U_\alpha$ and $[\alpha+1,\omega_1]\subseteq V_\alpha$.
Then $\gamma\N\setminus(U_\alpha\cup V_\alpha)$ is a compact subset of~$\N$,
and hence finite. By adding this finite set to $U_\alpha$ we can in fact
assume that $U_\alpha\cup V_\alpha=\gamma\N$ for all~$\alpha$.

We let $T_\alpha=U_\alpha\cap\N$ for all~$\alpha$.
The sequence $\omegaoneseq{T}$ has the following property:
\begin{quote}
$(*)$ if $\beta<\alpha$ then $T_\beta\setminus T_\alpha$ is finite 
and $T_\alpha\setminus T_\beta$ is infinite.
\end{quote}

Conversely, every such sequence determines a compactification of~$\N$.
To this end write $T_{\omega_1}=\N$ and then define a topology
on $\N\cup(\omega_1+1)$ as follows:
the points of~$\N$ are isolated and if $\alpha\le\omega_1$ then the sets
$$
(\beta,\alpha]\cup\bigl(T_\alpha\setminus(T_\beta\cup F)\bigr)
$$
where $\beta<\alpha$ and $F\subseteq\N$ is finite form a local base 
at~$\alpha$.
To ensure that the ordinal~$0$ is not an isolated point we tacitly assume
that $T_0$ is infinite.
The sets $U_\alpha=T_\alpha\cup[0,\alpha]$ and 
$V_\alpha=(\N\setminus T_\alpha)\cup[\alpha+1,\omega_1]$ are as in
the previous paragraph.

\bigskip
Before we investigate how the softness of~$\gamma\N$ depends on the
sequence $\omegaoneseq{T}$ we make a short digression on the cardinal
number~$\tee$ alluded to above.

A sequence $\omegaoneseq[\delta]{T}$ that satisfies property~$(*)$ above
is also called a tower.
The cardinal number~$\tee$ is defined to be the minimum ordinal~$\delta$
for which there is a \emph{complete tower}: a tower $\omegaoneseq[\delta]{T}$
with the additional property that if $S$~is such 
that $T_\alpha\setminus S$ is finite for all~$\alpha$ then 
$\N\setminus S$ is finite.

Thus, $\tee>\aleph_1$ means that \emph{every} tower
$\omegaoneseq{T}$ is \emph{incomplete}: there is an infinite set~$R$ such 
that $T_\alpha\cap R$ is finite, for all~$\alpha$, 
whereas $\tee=\aleph_1$ means that there is \emph{some} tower $\omegaoneseq{T}$
that is complete, that is, for which no such infinite set exists.

\bigskip
Now let $A$ and $B$ be disjoint subsets of~$\N$ such that 
$\cl A\cap\cl B\neq\emptyset$ in~$\gamma\N$.
We consider a few cases.

\smallskip
\noindent\textsc{Case 1}:
there is an~$\alpha\in\omega_1$ that is in the closure of $A$ and $B$.
Then $A$ and~$B$ have infinite intersections with $U_\alpha$ and
we can take infinite subsets~$C$ of~$A\cap U_\alpha$ and~$D$ 
of~$B\cap U_\alpha$ such that $\cl C=C\cup\{\alpha\}$ 
and $\cl D=D\cup\{\alpha\}$.  
Take any permutation~$h$ of~$\N$ that interchanges~$C$ and~$D$
and is the identity outside $C\cup D$.
Then $h$ extends to an autohomeomorphism of~$\gamma\N$ that is
the identity on~$\gamma\N\setminus\N$.
By construction $h[A]\cap B$ contains~$D$.
 
\smallskip
\noindent\textsc{Case 2}:
no such $\alpha$ can be found; hence $\cl A\cap\cl B=\{\omega_1\}$.
In this case we know that for every~$\alpha$ the intersections
$V_\alpha\cap A$ and $V_\alpha\cap B$ are both infinite.
Now we split into two subcases.

\noindent\textsc{Subcase 2a}:
there are infinite sets subsets~$C$ of~$A$ and~$D$ of~$B$ such that 
$C\cup D\subseteq^*V_\alpha$ for all~$\alpha$.
Then $C$ and $D$ converge to~$\omega_1$ and, as above, interchanging $C$ 
and~$D$ will witness softness.

By the definition of $\tee$ this subcase occurs, no matter how $\gamma\N$ is 
constructed, if $\tee>\aleph_1$.
This then proves the first statement at the beginning of this subsection.

\smallskip
\noindent\textsc{Subcase 2b}: 
one of $A$ and $B$ does not contain an infinite set as in \textsc{Subcase~2a}, 
say $A$ to be definite.
We claim that $\omega_1\cap\cl A$ is closed and unbounded in~$\omega_1$.
That it is closed is clear.
To show unboundedness let $\alpha\in\omega_1$.
The set $A\cap V_\alpha$ is infinite, but by assumption there is a $\beta$
such that $A\cap V_\alpha\setminus V_\beta$ is infinite.
It closure intersects~$\omega_1$ and is inside $V_\alpha\setminus V_\beta$,
hence $\cl A$ intersects the interal $[\alpha+1,\beta]$.

It now follows that $\omega_1\cap\cl B$ is closed and \emph{bounded} 
in~$\omega_1$; say $\omega_1\cap \cl B\subseteq[0,\alpha]$.
This means that $D=B\cap V_\alpha$ has $\omega_1$ as its only accumulation
point and hence that it converges to~$\omega_1$.

In this case the pair $\orpr DA$ witnesses non-softness: 
if $h:\gamma\N\to\gamma\N$ is a homeomorphism that is the identity 
on~$\omega_1$ then $h[D]$ converges to~$\omega_1$ and so 
$h[D]\subseteq^*V_\beta$ for all~$\beta$ and by our
assumption on~$A$ this means that $h[D]\cap A$ is finite.

Note that there is no \textsc{Subcase 2c}: in that subcase $\omega_1\cap\cl A$
and $\omega_1\cap\cl B$ would both be closed and unbounded and this would
bring us back to \textsc{Case 1}.

\smallskip
We shall show that under the assumption $\tee=\aleph_1$ one can construct
two compactifications, $\gamma_1\N$ and $\gamma_2\N$, of~$\N$ and with remainder
$\omega_1+1$ such that in~$\gamma_1\N$ \textsc{Case~1} always occurs,
so $\gamma_1\N$~is soft, and in $\gamma_2\N$ there is a pair~$\orpr AB$ where
\textsc{Subcase~2b} occurs, hence $\gamma_2\N$~is not soft. 

By the assumption $\tee=\aleph_1$ there is a complete tower~$\omegaoneseq{T}$. 

We let $\gamma_1\N$ be the compactification of~$\N$ determined by this
tower.
Let $A$ and $B$ be disjoint subsets of~$\N$ whose closures intersect
and assume $\omega_1\in\cl A\cap\cl B$.
Because neither $A$ nor $B$ contains an infinite subset as in 
\textsc{Subcase 2a} the argument in \textsc{Subcase 2b} applies 
to~$A$ and~$B$ to show that $\omega_1\cap\cl A$ and $\omega_1\cap\cl B$ 
are closed and unbounded in~$\omega_1$.
We find that we are in \textsc{Case 1} and that $\gamma_1\N$ is a soft
compactification of~$\N$.

To construct the promised non-soft compactification~$\gamma_2\N$ we take the 
sum $\gamma_1\N\oplus\alpha\N$ of~$\gamma_1\N$ and the one-point 
compactification~$\alpha\N$, and identify the points $\omega_1$ and $\infty$.

We let $A$ be the copy of~$\N$ from~$\gamma_1\N$ and
$B$ the copy of~$\N$ from~$\alpha\N$.
This pair witnesses \textsc{Subcase 2b} above and thus shows
that $\gamma_2\N$ is not a soft compactification of~$\N$.

\begin{remark}\label{rem.wt.aleph1}
As mentioned above the cardinal number $\tee$ is equal to the number~$\pee$
mentioned in the introduction.

A consequence of this is that we could have concluded right away that
$\tee>\aleph_1$ implies that every compactification of~$\N$ with
remainder~$\omega_1+1$ is soft.
We believe the argument that we gave above is more instructive.

Nevertheless we do record here for future use that under the assumption
$\tee>\aleph_1$ every compactification of~$\N$ whose remainder has 
weight~$\aleph_1$ or less is soft.
This means that if we want to prove, in~$\ZFC$, that some specific compact 
space of weight~$\aleph_1$ is soft-Parovichenko we can (or rather must)
assume that $\tee=\aleph_1$.
\end{remark}

\subsection*{The compact ordered space $\omega_1+1+\omega_1^\ast$}

The compact ordered space $K=\omega_1+1+\omega_1^\ast$ is of weight~$\aleph_1$
and hence is a Parovichenko space. 
It is only slightly more complicated than $\omega_1+1$ but, as we shall see,
its softness is already undecidable.

We think of~$K$ as the quotient of $(\omega_1+1)\times2$
obtained by identifying $\orpr{\omega_1}0$ and $\orpr{\omega_1}1$ to 
one point, which we call~$\Omega$.

Note that $(\omega_1+1)\times2$ is a soft remainder, as a sum of two soft
compactifications of~$\N$ is again a soft compactification.
This example will therefore show two things: 
a compact space of weight~$\aleph_1$ need not be a soft remainder, and
the continuous image of a soft remainder need not be a soft remainder itself.

It is easy to exhibit \emph{some} compactification of~$\N$ with remainder~$K$ 
that is not soft, assuming $\tee=\aleph_1$ of course: take the 
compactification~$\gamma_1\N$ from the previous subsection.
We take the quotient of $\gamma_1\N\times2$ obtained by identifying 
$\{\omega_1\}\times2$ to one point.
The sets $A=\N\times\{0\}$ and $B=\N\times\{1\}$ witness non-softness of
the resulting compactification: if $C\subseteq A$ is infinite then its 
closure contains at least one point of $\omega_1\times\{0\}$; any homeomorphism
that maps $C$ into~$B$ will move that point to $\omega_1\times\{1\}$.

We shall show that the principle (NT) from~\cite{MR1611327} implies that
every compactification of~$\N$ with remainder~$K$ contains two sets
that behave like $\N\times\{0\}$ and $\N\times\{1\}$ in the above example.

To formulate (NT) we need to define the notion of a weakly $\sigma$-bounded
family of infinite subsets of~$\N$: 
given a family~$\calA$ of infinite subsets of~$\N$ we 
let $\calA^\downarrow$ denote the family of infinite sets~$X$ for which there is 
a member of~$\calA$ that contains it.
We call $\calA$ \emph{weakly $\sigma$-bounded} if for every countable
subfamily~$\calX$ of~$\calA^\downarrow$ there is an~$A\in\calA$ such that
$A\cap X$ is infinite for all~$X\in\calX$.

The principle (NT) states the following:
\begin{quote}
for each weakly $\sigma$-bounded subfamily~$\calA$ of~$\pow(\N)$ 
and each subfamily~$\calB$ of~$\calA$ of cardinality at most~$\aleph_1$
there is a subset~$C$ of~$\N$ such that $C\cap B$ is infinite 
for all $B\in\calB$ \emph{and} for every infinite subset~$D$ of~$C$ there 
is an~$A\in\calA$ such that $A\cap D$ is infinite.
\end{quote}
In \cite{MR1611327} this principle was shown to be consistent with 
$\cee=\bee=\aleph_2$. 

\bigskip
Now let $\delta\N=\N\cup K$ be a compactification, where $K$~is the remainder.
For every~$\alpha$ we choose pairwise disjoint subsets 
sets $L_\alpha$, $M_\alpha$ and $R_\alpha$ of~$\N$ such that,
with closures taken in~$\gamma\N$
\begin{itemize}
\item $[0,\alpha]\times\{0\}\subseteq \cl L_\alpha$,
\item $[\alpha+1,\omega_1]\times2\subseteq \cl M_\alpha$, 
\item $[0,\alpha]\times\{1\}\subseteq \cl R_\alpha$, and
\item the three closures are pairwise disjoint.
\end{itemize}

We apply the principle~(NT) to the families
$\calL=\{L_\alpha:\alpha<\omega_1\}$ and
$\calR=\{R_\alpha:\alpha<\omega_1\}$,
and the subfamilies 
$\calB_L=\{L_{\alpha+1}\setminus L_\alpha:\alpha<\omega_1\}$ of~$\calL^\downarrow$
and $\calB_R=\{R_{\alpha+1}\setminus R_\alpha:\alpha<\omega_1\}$ 
of~$\calR^\downarrow$ respectively.

The families $\calL$ and $\calR$ are clearly weakly $\sigma$-bounded: 
if $\calB$ is a countable family of infinite sets such that for 
all $B\in\calB$ there is an~$\alpha_B$ with $B\subseteq L_{\alpha_B}$
then take $\alpha=\sup_B\alpha_B$; the set $L_\alpha$~is as required
because $B\subseteq^*L_\alpha$ for all $B\in\calB$.
The same argument works for~$\calR$.

The families $\calB_L$ and $\calB_R$ are of cardinality~$\aleph_1$
and refine $\calL$ and $\calR$, respectively.
The principle~(NT) then guarantees there are subsets $C_L$ and~$ C_R$
of~$\N$ such that
\begin{itemize}
\item $B\cap C_L$ is infinite, for all $B\in\calB_L$, and, likewise
      $B\cap C_R$ is infinite, for all $B\in\calB_R$, and
\item for every infinite subset~$D$ of $C_L$ (or $C_R$) there is 
      an $L\in\calL$ (or an $R\in\calR$) such that $D\cap L$ 
      (or $D\cap R$) is infinite
\end{itemize}
We derive some consequences from this.

For every $\alpha$ the set $L_{\alpha+1}\setminus L_\alpha$ converges to
the point $\orpr{\alpha+1}0$, hence $\orpr{\alpha+1}0\in\cl C_L$.
It follows that the point in the middle, $\Omega$, is in the closure
of~$C_L$. 
And by a symmetric argument $\Omega\in\cl C_R$ also.

The intersection $C_L\cap C_R$ is finite. 
For if it were infinite then by the second condition in~(NT) there
is an~$\alpha$ such that $L_\alpha\cap C_L\cap C_R$ is infinite,
and by a second application of that condition there is a~$\beta$
such that $L_\alpha\cap C_L\cap C_R\cap R_\beta$ is infinite.
But $L_\alpha\cap R_\beta$ is finite, contradiction.
So we may as well assume that $C_L$ and $C_R$ are disjoint. 

Now let $h$ be an autohomeomorphism of $\gamma\N$ with the property
that $h[C_L]\cap C_R$ is infinite.
Take an~$\alpha$ such that $h[C_L]\cap C_R\cap R_\alpha$ is infinite.
Then take a~$\beta$ such that 
$L_\beta\cap C_L\cap h\preim[C_R\cap R_\alpha]$ is infinite.
Take $\gamma\le\beta$ such that $\orpr\gamma0$~is in the closure of 
the latter set; then $h(\gamma,0)$ is in the closure of~$R_\alpha$.
Hence certainly $h(\gamma,0)\neq\orpr\gamma0$.

\subsection*{The Cantor cube $\Ccube$ and the Tychonoff cube $\Tcube$}

Two natural compact spaces of weight~$\aleph_1$ to consider are the
Cantor and Tychonoff cubes $\Ccube$ and~$\Tcube$.
We shall show that both are soft-Parovichenko.
Since the ordered space $\omega_1+1+\omega_1^\star$ is a subspace of
both cubes this shows that subspaces of soft-Parovichenko spaces
need not be soft-Parovichenko themselves.

As mentioned in Remark~\ref{rem.wt.aleph1} we can assume $\tee=\aleph_1$.

We apply Theorem~3 from~\cite{MR196693} and take a set $\calF$ of functions
from~$\N$ to~$\N$ that is of cardinality~$\cee$ and \emph{independent}, 
which means that given $f_1$, \dots, $f_k$ in~$\calF$ and $n_1$, \dots, $n_k$ 
in~$\N$ there is an~$m\in\N$ such that $f_i(m)=n_i$ for all~$i$.
One readily checks that this is equivalent to the following: the image
of the diagonal map $e:\N\to\N^\calF$ of the family~$\calF$ is dense
where we consider the product topology on~$\N^\calF$ and the discrete topology
on~$\N$.
 
We take an injective sequence $\omegaoneseq{f}$ of elements of~$\calF$
and a complete a complete tower $\omegaoneseq{T}$ in~$\N$.

We make a technical adjustment to these two sequences, as follows.

Let $N=\{\orpr mn\in\N^2: n\le m\}$. 
We define a new tower $\omegaoneseq{I}$ and a new independent sequence
$\omegaoneseq{F}$ of functions on~$N$.
For every $\alpha$ let 
\begin{itemize}
\item $I_\alpha=N\cap (T_\alpha\times\N)$, and
\item define $F_\alpha$ by
$$
F_\alpha(m,n)=
  \begin{cases}
  f_\alpha(n)&\text{if }m\notin T_\alpha\text{ and }f_\alpha(n)\neq0\text{, and}\\
  0         &\text{otherwise}   
  \end{cases}
$$
\end{itemize}
It is elementary to verify that $\omegaoneseq{I}$ is also a complete tower
and that the sequence $\omegaoneseq{F}$ is again independent.
We have additionally created some interplay between the two:
if $\beta\ge\alpha$ then $I_\alpha\setminus F_\beta\preim(0)$ is finite.

We adjust the two sequences once more, this time without renaming them.
We identify $N$ with~$\N$ via some bijection; and we let the codomains
of the functions~$F_\alpha$ be the set~$Q$ of rational numbers in~$[0,1]$,
via some bijection that sends~$0$ to~$0$.

We let $e:\N\to\Tcube$ be the diagonal map of the sequence $\omegaoneseq{F}$,
defined by $e(n)(\alpha)=F_\alpha(n)$.
By the remark above the image set~$e[\N]$ is dense in~$Q^{\omega_1}$ where
$Q$~carries the discrete topology, so $e[\N]$~is certainly dense in~$\Tcube$.

It follows that $\beta e:\betaN\to\Tcube$ induces a surjection from
$\betaNminN$ onto~$\Tcube$; this surjection determines a compactification 
of~$\N$ with $\Tcube$ as its remainder.

The compactification can be visualised as a subspace of the Alexandroff double
$A(\Tcube)$ of~$\Tcube$: the underlying set of $A(\Tcube)$ is $\Tcube\times2$
with $\Tcube\times\{0\}$ being the set of isolated points; 
our compactification~$\gamma\N$ of~$\N$ then is 
$$
(\Tcube\times\{1\})\cup (e[\N]\times\{0\}).
$$ 
We identify $\N$ with $e[\N]\times\{0\}$ and $\Tcube$ with $\Tcube\times\{1\}$.
We can save on notation by observing for $x\in\Tcube$ and $A\subseteq\N$ 
we have $x\in\cl A$ iff $x$~is an accumulation point of~$e[A]$;
and $A$~converges to~$x$ iff $e[A]$~converges to~$x$.

Let $A$ and $B$ be disjoint in~$\N$ and let $x\in\Tcube$ be in the intersection
of their closures in~$\gamma\N$.

Let $M$ be a countable elementary substructure of $H((2^{\aleph_1})^+)$ that
contains our tower, our independent sequence of functions, and $A$, $B$ and $x$.
Let $\delta=M\cap\omega_1$.

If $H$~is a finite subset of~$\delta$ and $\varepsilon>0$ 
then the basic open set 
$$
O(x,H,\varepsilon)=
\bigl\{y:(\forall\alpha\in H)\bigl(|y_\alpha-x_\alpha|<\varepsilon\bigr)\bigr\}
$$
meets $e[A]$ and $e[B]$ in infinite sets.
By elementarity there is an $\alpha\in\delta$ such that $I_\alpha$, 
and hence also $I_\delta$, has infinite intersections with these infinite sets.

Write $\delta$ as the union of an increasing sequence $\omegaseq{H}$ of finite
sets. 
By the remark above we can find sequences $\omegaseq{a}$ and $\omegaseq{b}$
in~$A$ and~$B$ such that 
$e(a_n)\in O(x,H_n,2^{-n})\cap A\cap I_\delta$
and $e(b_n)\in O(x,H_n,2^{-n})\cap B\cap I_\delta$ for all~$n$.

Let $y$ be the point defined by $y\restr\delta=x\restr\delta$ and 
$y(\alpha)=0$ for $\alpha\ge\delta$.
We claim that $\omegaseq{a}$ and $\omegaseq{b}$ converge to~$y$.
Indeed, let $H$ be a finite subset of~$\omega_1$ and~$\varepsilon>0$.
Fix $K$ such that $H\cap\delta\subseteq H_K$.
Also, using the fact that 
$I_\delta\setminus F_\alpha\preim(0)$ is finite for $\alpha\in H\setminus\delta$
find $L\ge K$ such that 
$F_\alpha(a_n)=F_\alpha(b_n)=0$ when $n\ge L$ and $\alpha\in H\setminus\delta$.
Then $n\ge L$ implies $a_n,b_n\in O(y,H,\varepsilon)$.

Now take the permutation $h$ of $D$ that interchanges $a_n$ and $b_n$ 
for all~$n$ and leaves the other elements in their places.

This shows that $\Tcube$ is soft-Parovichenko.

\bigskip
To show that $\Ccube$ is soft Parovichenko one can use basically the same 
argument as above.
The only change that needs to be made is to the diagonal map~$e$:
let
$$
e(n)(\alpha)=\begin{cases} 
              1 &\text{if } F_\alpha(n)\neq0\text{, and}\\
              0 &\text{if } F_\alpha(n)=0
              \end{cases}
$$

\section{Remarks and questions}

The of this paper uses the words `may be' rather than the word `are'
and the previous two sections show why that is.
Under the Continuum Hypothesis the `are' is justified but not in general.

In $\ZFC$ all compact spaces of weight~$\aleph_1$ are Parovichenko, but as 
we have seen the space $\omega_1+1+\omega_1^\star$ is a Parovichenko
space that is consistently not soft-Parovichenko.

This state of affairs suggests various further questions about the nature 
of soft remainders of~$\N$.

In exploring the possible parallels between the classes of Parovichenko
spaces and soft-Parovichenko spaces we saw that the Continuum Hypothesis 
simply implies that there is no difference.

The class of Parovichenko spaces is closed under continuous images;
the class of soft-Parovichenko spaces is not, consistently.

The class of Parovichenko spaces is, consistently, not closed under subspaces:
in the Cohen model where $\cee=\aleph_2$ the ordinal space~$\omega_2+1$
is not remainder of~$\N$ even though it is a subspace of the Parovichenko
space~$\Ccube[\omega_2]$.
We have seen that the same holds for soft-Parovichenko spaces.

A fair number of known Parovichenko spaces is also soft-Parovichenko;
see the list in the introduction.
Notably absent in that list are the separable compact spaces, so that will 
be our first question:

\begin{question}
Is every separable compact space soft-Parovichenko?  
\end{question}

There are a few spaces in this class that are worth singling out.
Given that every space of weight less than~$\tee$ is a soft remainder
we ask in particular

\begin{question}
Are the cubes $\Ccube[\tee]$ and $\Tcube[\tee]$ soft-Parovichenko?  
\end{question}

We can ask this for every cardinal in the interval $[\tee,\cee]$ but instead
we ask whether there is some relationship between these cardinals.

\begin{question}
If $\kappa<\lambda$ and $\Ccube[\lambda]$ is a soft remainder is then
$\Ccube[\kappa]$ soft as well?
Likewise for the Tychonoff cubes.  
\end{question}

\end{document}